\documentclass[sn-aps]{sn-jnl}
\theoremstyle{thmstyleone}%
\usepackage{booktabs,algpseudocode}
\usepackage{amssymb}
\usepackage{amsmath}
\usepackage{graphicx}
\usepackage{cleveref}
\theoremstyle{thmstyletwo}%
\theoremstyle{thmstylethree}%
\begin{document}
\title[Hierarchical deep learning based adaptive time stepping]{Hierarchical deep learning based adaptive time stepping scheme for multiscale simulations}

\author[1]{\fnm{Asif} \sur{Hamid}}
\author*[1]{\fnm{Danish} \sur{Rafiq}}\email{danish.rafiq@iust.ac.in}
\author[1]{\fnm{Shahkar Ahmad} \sur{Nahvi}}
\author[2]{\fnm{Mohammad Abid} \sur{Bazaz}}
\affil[1]{\orgdiv{Department of Electrical Engineering}, \orgname{Islamic University of Science and Technology}, \orgaddress{\street{Awantipora}, \postcode{192122}, \state{Jammu \& Kashmir}, \country{India}}}
\affil[2]{\orgdiv{Department of Electrical Engineering}, \orgname{National Institute of Technology,} \orgaddress{ \city{Srinagar}, \postcode{190006}, \state{Jammu \& Kashmir}, \country{India}}}

\abstract{Multiscale is a hallmark feature of complex nonlinear systems. While the simulation using the classical numerical methods is restricted by the local \textit{Taylor} series constraints, the multiscale techniques are often limited by finding heuristic closures. This study proposes a new method for simulating multiscale problems using deep neural networks. By leveraging the hierarchical learning of neural network time steppers, the method adapts time steps to approximate dynamical system flow maps across timescales. This approach achieves state-of-the-art performance in less computational time compared to fixed-step neural network solvers. The proposed method is demonstrated on several nonlinear dynamical systems and source codes are provided for implementation. This method has the potential to benefit multiscale analysis of complex systems and encourage further investigation in this area.}


\maketitle
\vspace{10cm}
\section{Introduction}
\label{sec1}
\par
Multiscale systems are ubiquitous in science and engineering. Modeling and controlling such systems is essential due to their prevalence in natural and engineered systems, and understanding their behavior requires a multidisciplinary approach that integrates models, and experimental techniques at multiple scales \cite{weinan2011book}. These complex systems generally have dynamics operating at different spatiotemporal scales, such as a \textit{fine} or microscale, and a \textit{coarse} or macroscale. Microscale modeling usually involves analyzing the system behavior at fine resolutions, thus entailing a substantial computational cost while capturing the system's long-term behavior. On the other hand, the macroscopic models are efficient, but their accuracy hinges on the ability to capture the system dynamics effectively. Another challenge in studying multiscale systems is that the governing equations may be explicitly known at the microscopic/individual level, but the closures required to translate them to high-level macroscopic descriptions remain elusive. For instance, at the microscale, the governing equations that describe the behavior of fluid particles can be modeled using molecular dynamics simulations, which consider the interactions between individual particles. However, at the macroscale, the behavior of fluid flows is mainly influenced by interactions between particles at larger length scales, such as turbulent eddies, which are challenging to model accurately due to their complex nature. Thus, multiscale analysis involves deriving representative models at different scales and coupling them to achieve the accuracy of the microscopic models, as well as the efficiency of the macroscopic models \cite{weinan2011book}.
\par 
 Many efforts have been made towards this goal, and various multiscale modeling techniques have been developed that combine different approaches to analyze and study the behavior of such systems. Some of these classical methods include the equation-free method (EFM) \cite{kevrekidis2003equation,kevrekidis2004equation}, multi-grid methods \cite{mccormick1987multigrid}, the heterogeneous multiscale method (HMM) \cite{weinan2003heterognous}, and the flow averaged integrator (FLAVOR) method \cite{tao2010nonintrusive}. While EFM and HMM use \textit{coarse time steppers} to simulate the evolution of macroscopic variables through a microscopic simulation, the FLAVOR method uses averaging flows to study the system behavior. However, the accuracy of these methods is highly dependent on the separation between different scales, the type of time integrator used, and how well the information is captured across different scales. On the other hand, data-driven modeling techniques, such as those based on the \textit{Koopman} theory, have also gained much attention to model multiscale phenomena, mainly when the governing equations or closure models are unavailable or difficult to derive. Most of these methods aim to separate the complex multiscale time-series data into its constituent timescale components. These include the multi-resolution dynamic mode decomposition (DMD) method \cite{mrdmd}, the sliding DMD-based method \cite{dylewsky2019dynamic}, sparse DMD method \cite{manohar2019optimized}, the transfer-operator-based methods \cite{transfer1,transfer2}, and the system identification-based method \cite{champion2019discovery}. These techniques primarily use a windowed subsets of the data to recover local linear models, and the dominant timescales are then identified via the spectral clustering of the eigenvalues. These methods however encounter several challenges, including problems in handling noisy data, inefficient convergence outside of training data, and the curse of dimensionality for large-scale systems \cite{baddoo2023physics}. 
 \par 
 Another recent approach to handle complex multiscale systems is based on machine learning (ML) techniques \cite{wehmeyer2018time,liu2023multiresolution,raissi2019physics,alber2019integrating,vinuesa2022enhancing}. These methods typically perform well due to the remarkable performance of deep neural networks (DNNs). As neural networks are universal approximators, these are used to approximate any continuous function with sufficiently many hidden units \cite{weinan2017proposal}. For instance, ML methods have been used to model dynamical systems \cite{milano2002neural,zhang2021midphynet,bailer1998recurrent}, in reduced modeling approaches \cite{lee2020model,vlachas2022LED,nakamura2021convolutional}, for attractor reconstructions \cite{pathak2017using,lu2018attractor}, in forecasting applications \cite{vlachas2018data,wiewel2019latent,zhang1998forecasting}, in mesoscopic material modeling \cite{fish2021mesoscopic}, in biological systems \cite{meier2009multiscale,weinberg2010multiscale}, and in molecular kinematics \cite{wehmeyer2018time,mardt2018vampnets,vlachas2021accelerated}. ML-based methods have also been merged with classical numerical techniques to perform discrete-time stepping \cite{parish2020time,regazzoni2019machine,qin2019data,raissi2018multistep,rudy2019deep,kim2020robust}.  With regard to multiscale systems, the flow map viewpoint of dynamical systems \cite{qin2019data,ying2006phase} has been used for model discovery of multiscale physics \cite{bramburger2020poincare,bramburger2020sparse}. This idea has been further exploited to build the multiscale hierarchical time-stepping (HiTS) method for multiscale problems in Ref. \cite{liu2022hierarchical}. In particular, the authors use a hierarchy of neural network time stepper (NNTS) models to approximate the dynamical system flow map over a range of timescales. This has the advantage that the NNTS models are unconstrained by the local \textit{Taylor} series expansion, unlike the classical time-steppers, and can capture a range of timescales with high accuracy. However, this method involves a computationally expensive cross-validation step for shortlisting the NNTS models for prediction and employs a fixed step size-based time-stepping strategy which is inefficient due to numerical consideration.
\par
In this contribution, we build upon the multiscale HiTS method and propose an adaptive multiscale machine learning framework that balances computational efficiency and accuracy. Notably, we demonstrate how a hierarchy of NNTS models, trained at different timescales, can be adaptively used during the simulation based on the multiscale properties of the system. This has the advantage that the proposed adaptive hierarchical time-stepping (AHiTS) strategy achieves the same accuracy as multiscale HiTS \cite{liu2022hierarchical} but with fewer time steps entailing computational savings. Also, the prediction accuracy can be adjusted by automatically adjusting the NNTS models on-the-fly during the simulation. Our work is motivated by the variable time-step solvers \cite{atkinson2011numerical} that adjust the relative time-step based on a given error tolerance. In summary, the main contributions of this work are:
\begin{enumerate}
	\item We present a novel deep learning-based adaptive time-stepping scheme for multiscale systems that balances accuracy and efficiency (in Section \ref{ahits}).
	\item The proposed AHiTS scheme uses a hierarchical deep learning perspective to train neural network time stepper models at different timescales. These models are then used adaptively as per the multiscale dynamics of the system in an efficient manner.
	\item The resulting AHiTS framework achieves the same level of accuracy as current state-of-the-art multiscale HiTS method \cite{liu2022hierarchical} in less computational time.
	\item Each NNTS model is trained independently and can focus on a given timescale for a short period, thus avoiding the exploding/vanishing gradient problem in neural networks.
	\item The AHiTS framework is robust to noise and works in a data-driven framework compared to classical physics-based multiscale methods.
	\item We also provide the necessary source codes for training the networks in an open-source Python environment.
\end{enumerate}
The remainder of this paper is sectioned as follows: In Section \ref{ts}, we provide a tutorial overview of hierarchical deep learning method for multiscale systems. We also describe how neural network-based flow map learning works similar to an explicit \textit{Euler} time-stepping scheme. Then, in Section \ref{ahits}, we present the proposed AHiTS scheme based on an adaptive selection of NNTS models. We give the flow chart of the proposed method and provide numerical algorithms to train and test the networks. Next, in Section \ref{sec-sims}, we demonstrate the application of AHiTS on several canonical ODEs and PDEs. We thoroughly discuss the results and provide a comparative analysis with the multiscale HiTS scheme \cite{liu2022hierarchical}. Finally, in Section \ref{sec-discuss}, we provide a detailed discussion on various aspects of the proposed method including the case of noisy experiments and the few limitations of this approach.
\section{Hierarchical deep learning using residual neural networks}
\label{ts}
In this section, we provide an overview of the multiscale HiTS method proposed in Ref. \cite{liu2022hierarchical}. In the next section, we build on this and demonstrate how adaptive selection of different NNTS models can yield an efficient simulation scheme for multiscale problems. To proceed with the idea, consider a nonlinear, continuous-time, multiscale dynamical system of the form:
\begin{equation}
	\dot{\mathbf{x}}(t)= \mathbf{f}(\mathbf{x}(t),t), \quad \mathbf{x}(t_0)=\mathbf{x}_0,
	\label{sys}
\end{equation}
where $\mathbf{x}(t) \in \mathbb{R}^n  $ represents the state vector of the given system with $n$ degrees of freedom, $t \in [0, ~ t_f]$ denotes the time, $\mathbf{x}_0$ is the initial condition of the state vector $\mathbf{x}(t)$, and $ \mathbf{f}(:,t) : \mathbb{R}^{n} \rightarrow \mathbb{R}^{n}$  represents the nonlinear mapping describing the evolution of the system. In most of the engineering applications, the ordinary differential equation (ODE) system (\ref{sys}) usually results from the spatial discretization of the underlying partial differential equation (PDE) via the finite difference method (FEM), finite element method (FEM) or using the spectral methods \cite{conte2017elementary}. Since most of these systems are inherently nonlinear, the analytical solution are rarely available, as such, the solution of (\ref{sys}) is numerically obtained by employing different time-stepping integration schemes such as explicit/implicit, fixed/adaptive, and one-step/multistep methods with varying degrees of stability \cite{atkinson2011numerical}. These numerical methods take successive time steps in an iterative manner to construct the solution trajectory. For instance, given the initial condition $\mathbf{x}(t_{0})$, these methods approximate the discrete-time flow map
\begin{equation}
	\mathbf{x}(t+\Delta t) =\mathcal{N}(\mathbf{x}(t),\Delta t)\triangleq \int_{\tau=t}^{t + \Delta t} \mathbf{f}(\mathbf{x}(\tau),\tau) d\tau,
	\label{int}
\end{equation}
usually via the \textit{Taylor-series} expansion \cite{butcher1987numerical,guckenheimer2013nonlinear,wiggins2003introduction}. However, the accuracy of these methods depends upon the local step size $\Delta t$ employed during the integration making such time discretization
local in nature. For instance, the classical \textit{ode45} method uses a fixed-step, fourth-order \textit{Runga-Kutta} scheme that incurs a local truncation error of $\mathcal{O}(\Delta t^{5})$ at every step of the solution, and a global truncation error of $\mathcal{O}(\Delta t^{4})$. In contrast to this, neural networks-based time-steppers remain unconstrained by the step size constraints as they rather learn the underlying discrete flow map that progresses the states forward in time \cite{qin2019data}. Among various neural network architectures, residual neural networks (ResNets) have been widely used because of the lesser training times and have been identified to work similarly to a fixed-step, first-order \textit{Euler's} scheme (see e.g., Refs. \cite{chang2017multi,qin2019data,chen2018neural,he2016deep}) and are discussed as follows.

\subsection{Flow map learning via ResNets}
\begin{figure}[htbp]
	\centering
	\includegraphics[width=0.6\textwidth]{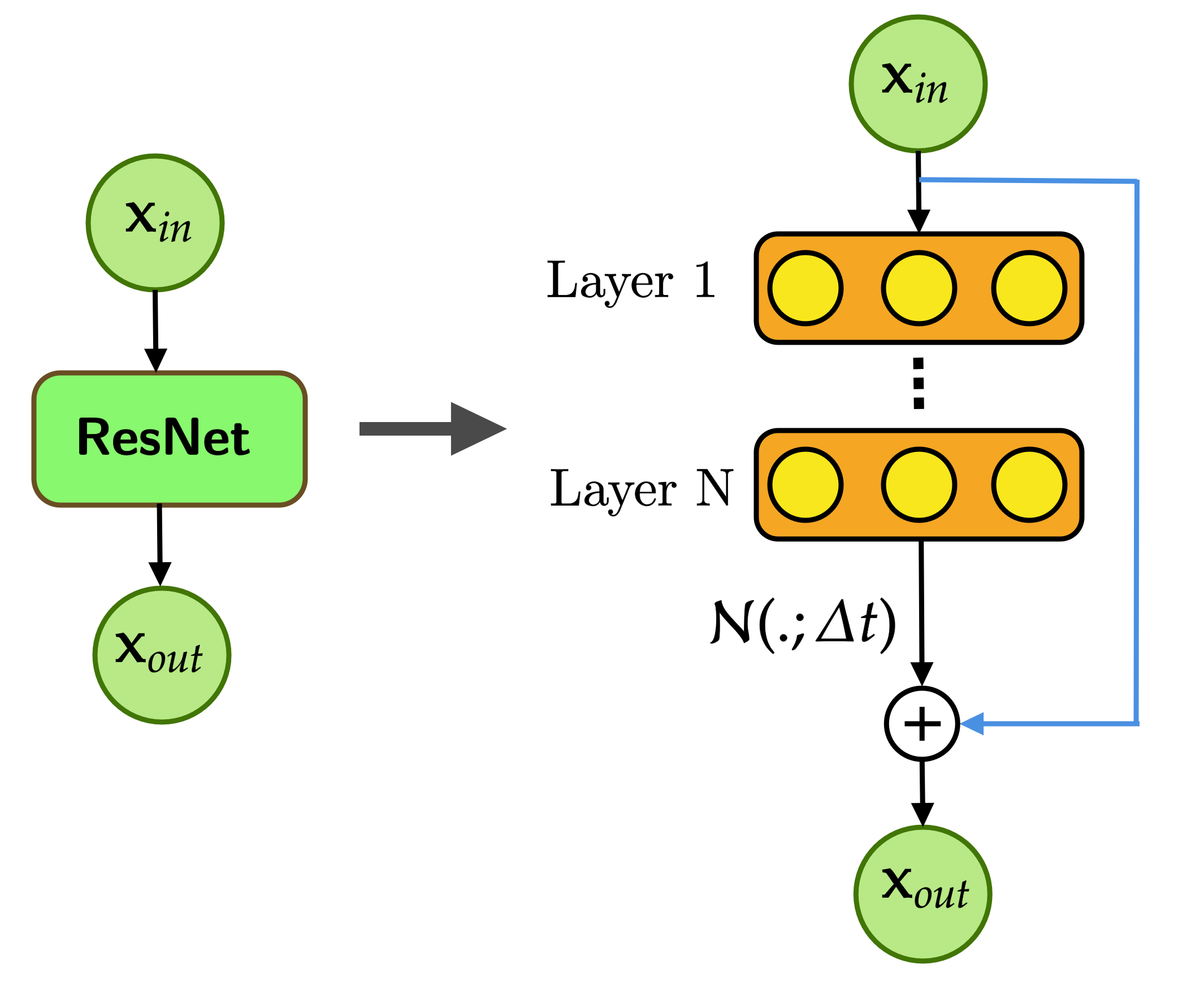}
	\caption{Schematic for one step approximation using ResNet}
	\label{resnet}
\end{figure}

Residual neural networks (ResNets) are a type of neural network architecture that aim to solve the problem of \textit{vanishing gradients} in DNNs \cite{goodfellow2016deep}. This issue arises when the gradient signal becomes too small during backpropagation, making it difficult for the network to learn and adjust its parameters effectively. ResNets address this problem by introducing skip connections, which allow the network to learn residual functions that make it easier to propagate gradients. To learn the flow map of system (\ref{sys}) using ResNets, we proceed as follows. We begin by collecting $p$ trajectories sampled at $t_f$ instances with time step $ \Delta t$ given as:
\begin{equation}
	\mathbf{X}^{i}=\Bigg[ \mathbf{x}^{i}(t_0)~ \quad \mathbf{x}^{i}(t_0+ \Delta t)~ \dots ~ \mathbf{x}^{i}({t_0}+t_f\Delta t)~ \Bigg], 
	\label{data}
\end{equation}
where $ i \in 1,\dots,p$. This data set is then used to train a ResNet model by feeding the first $t_{f}-1$ entries of $\mathbf{X}$, i.e., $[\mathbf{x}^{i}(t_0), ..., \mathbf{x}^{i}({t_0}+(t_{f}-1)\Delta t)]$ as the input and one $\Delta t$ ahead entries, i.e., $[\mathbf{x}^{i}(t_0+\Delta t), ... ,\mathbf{x}^{i}({t_0}+t_f\Delta t)]$ as the output. By doing so, the network learns the $\Delta t$ flow map for system (\ref{sys}) given as:
\begin{equation}
	\mathcal{N}(.;\Delta t)=\mathbf{a}_L(\mathbf{W}_{L}(\dots \mathbf{a}_1(\mathbf{W}_1))\dots)),
	\label{NNF}
\end{equation}
where $L$ is the number of layers, $\mathbf{a}_{k}~ (k=1,..,L)$ is the type of activation functions used at layer $k$, and $\mathbf{W}_{k}$ corresponds to weights for $k^{th}$ layer. The corresponding weights are learned by minimizing the loss function
\begin{equation}
	\mathbf{W}^{*}= \arg \min_{\mathbf {W}} \frac{1}{p\times t_f} \sum_{i=1}^{p} \sum_{j=1}^{t_f}  \left[\mathcal{L}( \mathbf{\hat{x}}^{i}_{t+j\Delta t},\mathbf{x}^{i}_{t+ j \Delta t})
	\right],
	\label{opt}
\end{equation}
where $ \mathcal{L}$ denotes the discrepancy error between the true state $ \mathbf{x}$ and the predicted state $\mathbf{\hat{x}}$. Figure \ref{resnet} shows the pictorial diagram of flow map leaning via a ResNet. As can be seen, the ResNet block is composed of $L$ hidden layers, and an identity operator that adds the current input $\mathbf{x}_{in}$ to its output $\mathbf{x}_{out}$ given as:
\begin{equation}
	\mathbf{x}_{out}=\mathbf{x}_{in} + \mathcal{N}(.;\Delta t),
	\label{resneteq}
\end{equation}
where $\mathbf{x}_{in}=\mathbf{x}(t_{0})$ and $\mathbf{x}_{out}=\mathbf{x}(t_{0}+\Delta t)$. The above equation (\ref{resneteq}) resembles to that of the forward \textit{Euler's} discretization method that queries the nonlinear function $\mathbf{f}(\mathbf{x}(t),t)$ in system (\ref{sys}) at every time step $\Delta t$ to obtain the next step solution. Thus, similar to a forward time-stepping algorithm that approximates the discrete-time flow map of a dynamical system, a neural network time stepper is able to learn this underlying flow map and provide predictions at multiple steps forward in time for unknown initial states. The advantage with neural networks is that they work in an equation-free manner and remain unconstrained by the stability margins of the local step sizes.
\subsection{Coupling multiple ResNets for multiscale systems}
\label{sec-couple}
\begin{figure}[H]
	\centering
	\includegraphics[width=1\textwidth]{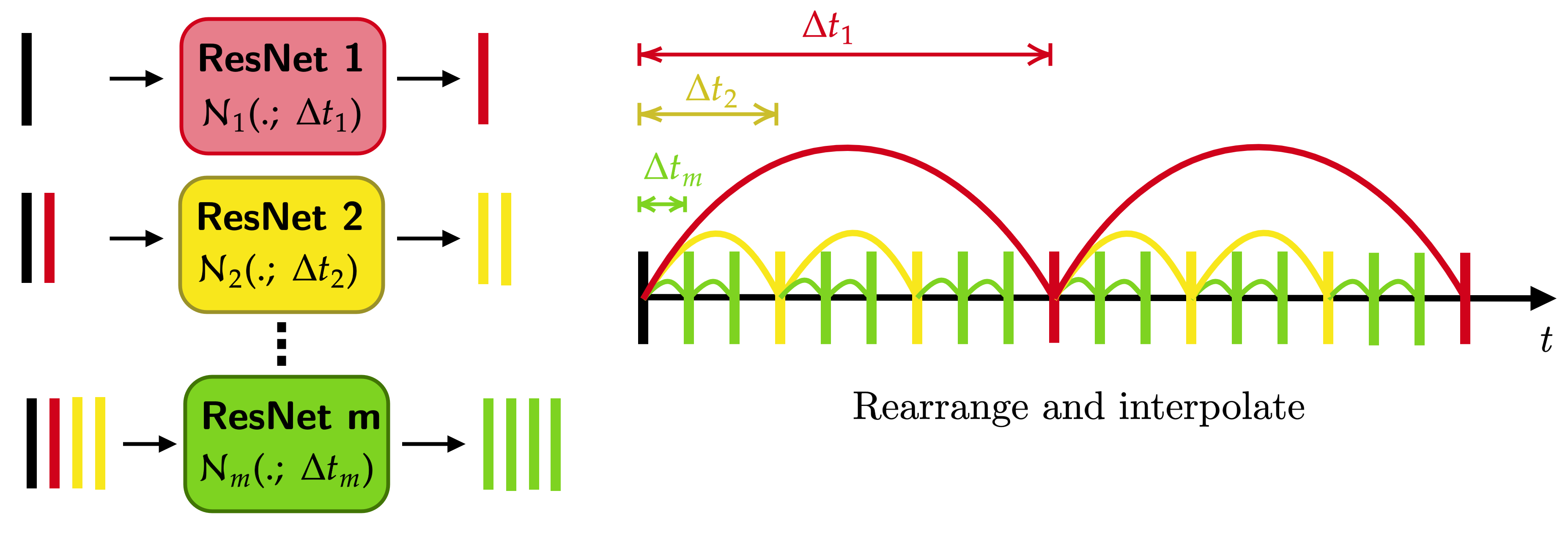}
	\caption{Coupling various ResNets for multiscale simulations (adapted from Ref. \cite{liu2022hierarchical})}
	\label{resnetmul}
\end{figure}

As suggested in Ref. \cite{liu2022hierarchical}, different ResNets models trained at different timescales can be coupled together for multiscale applications, i.e., the NNTS models trained at various levels can be used independently to evaluate the fast and slow dynamics of the systems. The idea is illustrated in Fig. \ref{resnetmul} which shows the basic schematic for multiscale HiTS method. As can be seen, $m$ different NNTS models i.e., $\mathcal{N}_{d}(.;\Delta t_d )~ (d=1,..,m)$ are trained at multiple timescales of the unit time-step $\Delta t$, i.e., $\Delta t_{d}= 2^{m-d} \Delta t$. This results in capturing the long-term behavior with the NNTS model having larger time-step and short-term behavior with smaller ones. Another advantage is that the training can be achieved independently to make the computations inherently parallelizable. However, before coupling the individual responses, a cross-validation step is performed to filter the NNTS models that yield the best time-stepping performance, i.e., given $m$ neural network models $\{\mathcal{N}_{1},..,\mathcal{N}_{m}\}$, ordered by associated step sizes, the models are filtered by imposing a lower bound $u$ and an upper bound $v$ $\{\mathcal{N}_{u},..,\mathcal{N}_{v}\}$
for all $m$. Once the models are selected, vectorized computation is performed to record the predictions from the models (see Algorithm 1 from Ref. \cite{liu2022hierarchical}).
 Coupling various NNTS models for multiscale modeling is a fairly decent strategy; however, it has some underlying limitations. Firstly, cross-validating different models before coupling incur huge computational costs, especially for large-scale models. This is because of the exhaustive search strategy that checks every possible combination of the models. Secondly, the direct coupling method results in a fixed step size-based prediction, i.e., the smallest and largest steps-sizes are set by the lower and upper bound $u$ and $v$ on the models. This has the limitation that time steps cannot be adaptively changed during the simulation resulting in an inefficient computation. Motivated by these shortcomings, we propose an adaptive time-stepping scheme based on NNTS models which is discussed next.

\section{Deep learning based adaptive time stepping scheme}
\label{ahits}
\begin{figure}[t]
	\centering
	 \includegraphics[width=1\textwidth]{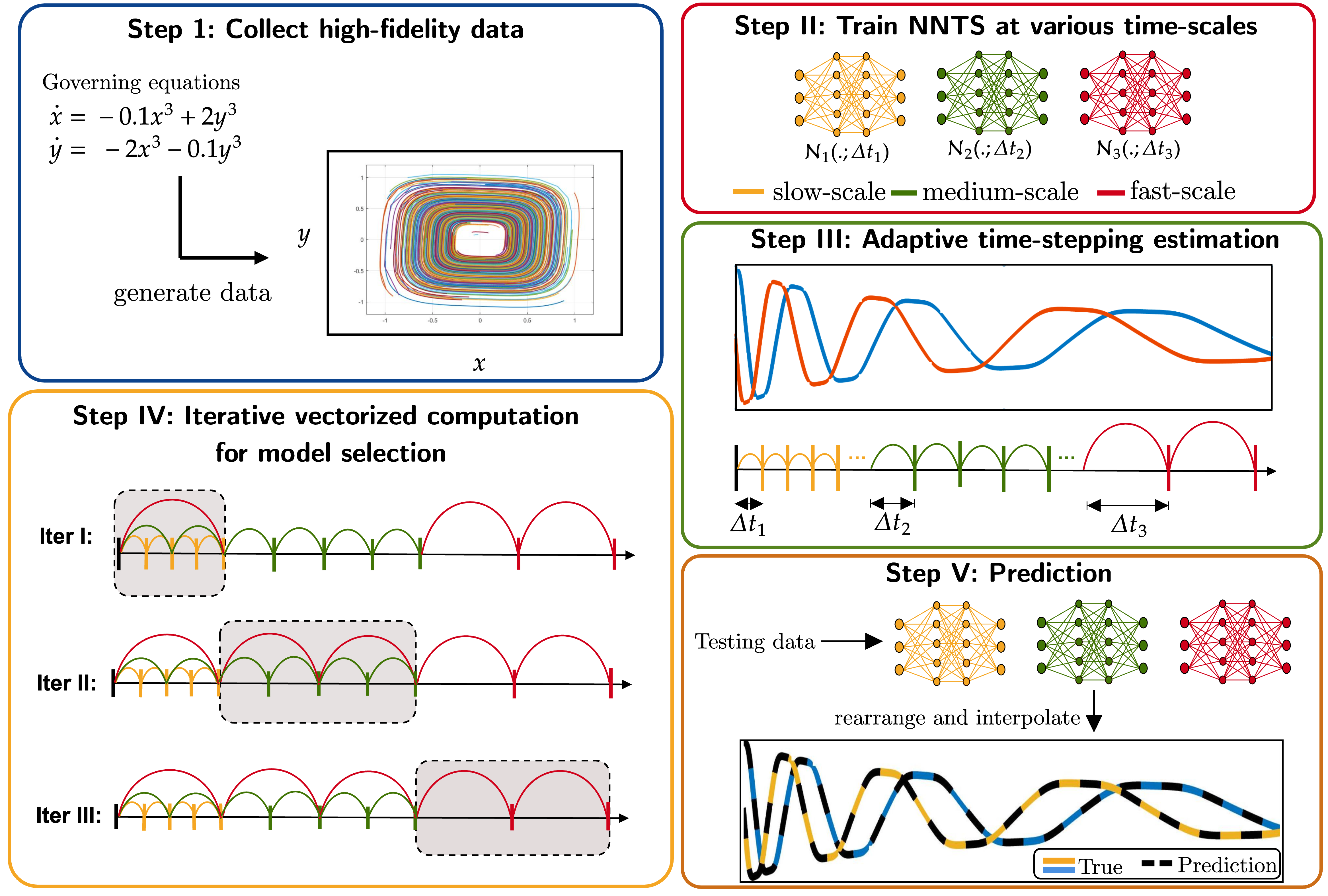}
	\caption{Block diagram of AHiTS algorithm: (Step I) Collection of data at different step sizes, (Step II) Training NNTS, (Step III) Calculation of step size according to system dynamics, (Step IV) Piece-wise vectorisation, (Step V) Predicting the output response on random test data.}
	\label{ahits_block}
\end{figure}
Here we outline the proposed deep-learning based adaptive time-stepping method. The overall schematic is presented in Fig. \ref{ahits_block}. The proposed AHiTS method has five main steps, which are explained using a simple cubic oscillator example.
\subsection{Data collection} As discussed in the previous section, we start by collecting the measurements of the system for a specified domain of interest in the state space with uniformly varying initial states of training, validation, and testing trajectories. Figure \ref{ahits_block} shows the phase portrait of the two-dimensional cubic oscillator system. The training data collected at unit time-step $\Delta t$ is sub-sampled uniformly at different multiples of $\Delta t$, i.e., $\Delta t_{d}=2^{m-d}\Delta t, d=1,2,..,m$ to train the individual ResNets. The validation data set is used to cross-validate the necessary time steps for automation, and the testing data is used as a reference to compare the performance of the method. We also add measurement noise to the data to test the method's robustness to noise before training the networks, which is discussed in Section \ref{sec-noise}. 
\subsection{Training NNTS models at different timescales} This step is similar to the multiscale HiTS method described previously in Section \ref{sec-couple}. The data collected in the previous step is used to train the individual ResNets for learning multiple flow maps $\mathcal{N}_{d}(.;\Delta t_{d})$ to capture the system dynamics across timescales. Kindly note that the choice of the hyper-parameters, such as the width and depth of the networks, play a crucial role in learning flow maps better. As suggested in \cite{liu2022hierarchical}, both deeper and wider networks allow capturing complex flow maps, including maps with larger temporal gaps. The training error also provides a good indication of how well a neural network represents a particular flow map.
\subsection{Adaptive time-stepping estimation}
Once the models are trained, the adaptive time-stepping is achieved as follows. Given an initial state $\mathbf{x}(t_{0})$, the NNTS model with largest $\Delta t$ is used to derive the next prediction at $\mathbf{x}(t_{0}+2^{m}\Delta t)$. The state vector $\mathbf{x}(t)$ evolution is tracked and compared against a predetermined tolerance $\epsilon$. If the evolution of the states is slow, i.e., the mean square error of $\|\mathbf{x}(t_{0}+2^{m}\Delta t) - \mathbf{x}(t_{0}) \| < \epsilon$, the step is recorded, and the same NNTS will be used to calculate the next prediction at $\mathbf{x}(t_{0}+2^{m+1}\Delta t)$. On the other hand, if the states evolve rapidly, this step is discarded, and the next NNTS model with a smaller time step will be checked. The process is repeated until an appropriate NNTS model with a reasonable step size is obtained. This procedure results in a variable time-step procedure where an NNTS model with a larger temporal gap captures the slow-evolving dynamics and vice-versa. The overall process is repeated until the current step reaches the final time step. The flowchart for this procedure is shown in Fig. \ref{ahits_flow}, and the corresponding algorithm is presented in Algorithm \ref{alglo1} of Appendix \ref{sec-algo}. Figure \ref{ahits_block} illustrates the adaptive time-stepping procedure for the case of the cubic oscillator, wherein smaller time steps are initially taken due to the fast-evolving dynamics and then gradually increased as the system gets slower. We can observe that this information is beautifully captured by adaptively refining the NNTS models. The choice selection of tolerance $\epsilon$ will be discussed in Section \ref{sec-discuss}. 
\subsection{Iterative vectorized computations}
One drawback with adaptively changing the neural network time steppers and the associated time steps is that when the system evolves very rapidly, the network with a smaller time step is used continuously. Although networks will smaller time step offer accurate short-term predictions, the error accumulates at every step and quickly dominates. To avoid error propagation, we employ the vectorized computations proposed in Ref. \cite{liu2022hierarchical}. The basic idea is to start by using the NNTS model with the highest temporal gap, like in the previous step. The state predictions obtained via this network are then stacked with original states and fed to the next-level NNTS model to initialize its hidden states. This is repeated for all the time-steppers. However, in our case, we iteratively perform this task to preserve the steps obtained in the previous step. In the first iteration, a time window selects NNTS models with similar time steps, and then vectorized computations are performed to record predictions. In the next iteration, this window selects the next group of time steppers, and the process repeats for all the networks. This strategy prevents local error accumulation while preserving the time steps obtained at the previous stage. Figure \ref{ahits_block} shows the schematic for the cubic systems where a windowed vectorization in a moving window is performed. The flowchart for performing the iterative vectorized computation is presented in Fig. \ref{ahits_flow} and the corresponding algorithm in Algorithm \ref{alglo2} of Appendix \ref{sec-algo}. It is worth noting that if the time window is kept equal to the entire time duration, and a high tolerance is used, then the proposed AHiTS coincides with the multiscale HiTS scheme, i.e., the networks from the AHiTS method are the same as shortlisted by the multiscale HiTS method. Thus, the proposed scheme can be viewed as a generalized multiscale simulation strategy.
\subsection{Prediction}
The final stage involves testing the adaptive strategy for unseen initial conditions. The predictions are recorded at time steps obtained during the adaptive time stepping estimation followed by the iterative-vectorization scheme, as explained above. Finally, a linear interpolation method is used to estimate the states at time steps that are not directly obtained from the AHiTS scheme.
\par In the next section, we will demonstrate the application of the AHiTS scheme on various benchmark models and provide a comprehensive analysis with the multiscale HiTS method.
\begin{figure}[H]
	\centering
	\includegraphics[width=\textwidth]{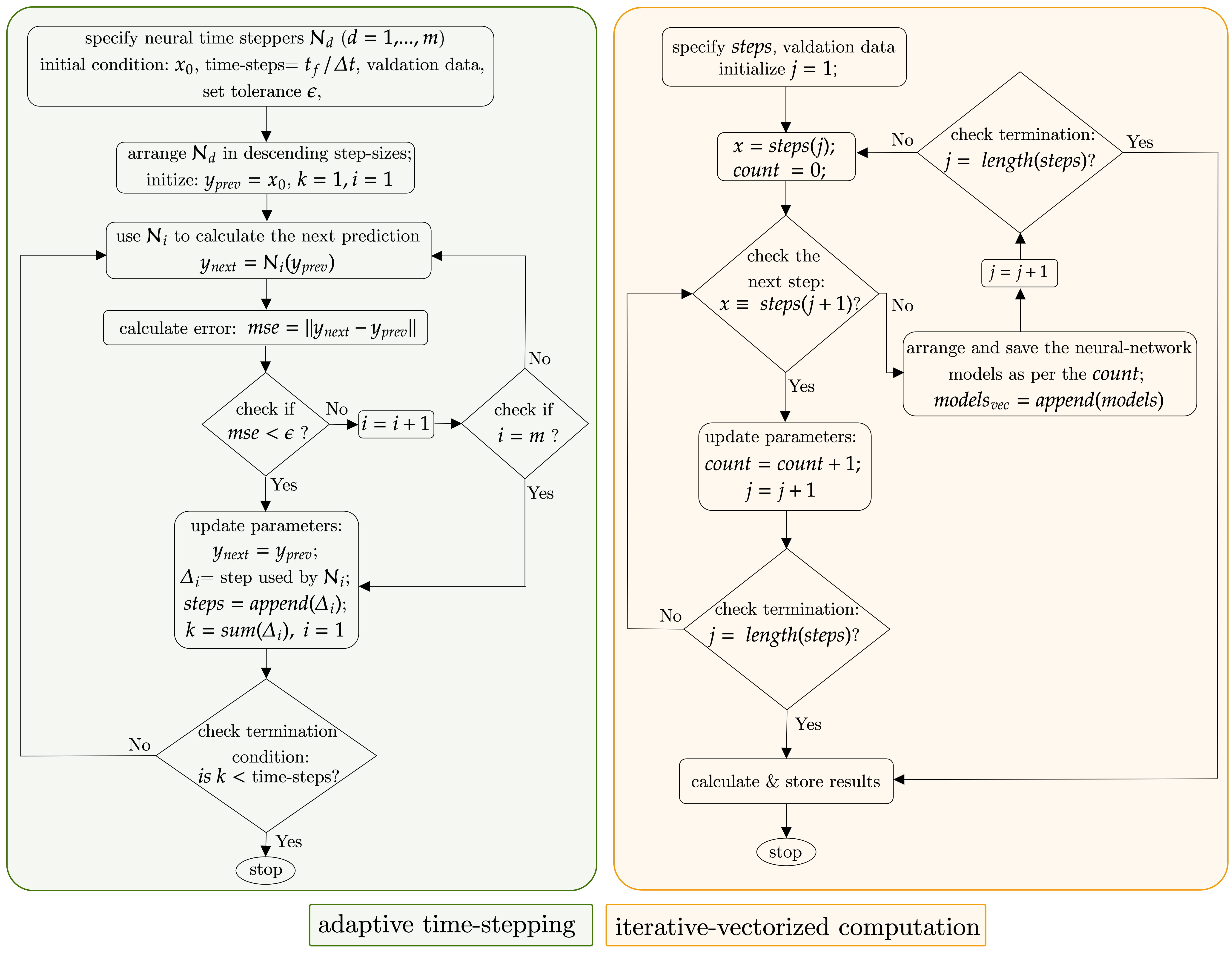}
	\caption{Flowcharts for steps 3 and 4 of the proposed AHiTS scheme}
	\label{ahits_flow}
\end{figure}

\section{Numerical Simulations}
\label{sec-sims}
In this section, we present the numerical simulations to demonstrate the advantages of the proposed AHiTS algorithm. To facilitate comparison with the multiscale HiTS method \cite{liu2022hierarchical}, all the model parameters and data structures are kept the same. The unit time step is fixed to be $\Delta t=0.01s$ for all the cases, and for testing the accuracy of the methods, we use the mean squared error (MSE) as a metric for comparison. All the simulations are performed using the open source Python API for PyTorch framework running on hp workstation Z1 with $11^{th}$ Gen Intel(R) Core (TM) $i7$-$11700$ $@2.5$GHz CPU. However, the code is written to take the advantage of faster computations via GPU if available.
\subsection{Benchmark ODEs}
We first benchmark the adaptive hierarchical time stepping method against four canonical nonlinear systems: a nonlinear system with a hyperbolic fixed point, a two-dimensional damped cubic oscillator, the Van der Pol oscillator, and a Hopf normal form. For all cases, we train eleven individual neural network time steppers, i.e., NNTS 0 - NNTS 10 with corresponding step sizes $\Delta t_{d}=2^{d}\Delta t ~(d=0,..,10)$ and couple them in our AHiTS scheme. We also combine the models as per the multiscale HiTS method. All the network parameters used for these systems are reported in Appendix Table \ref{parode}.
\begin{figure}[H]
	\centering
	\includegraphics[width=\textwidth]{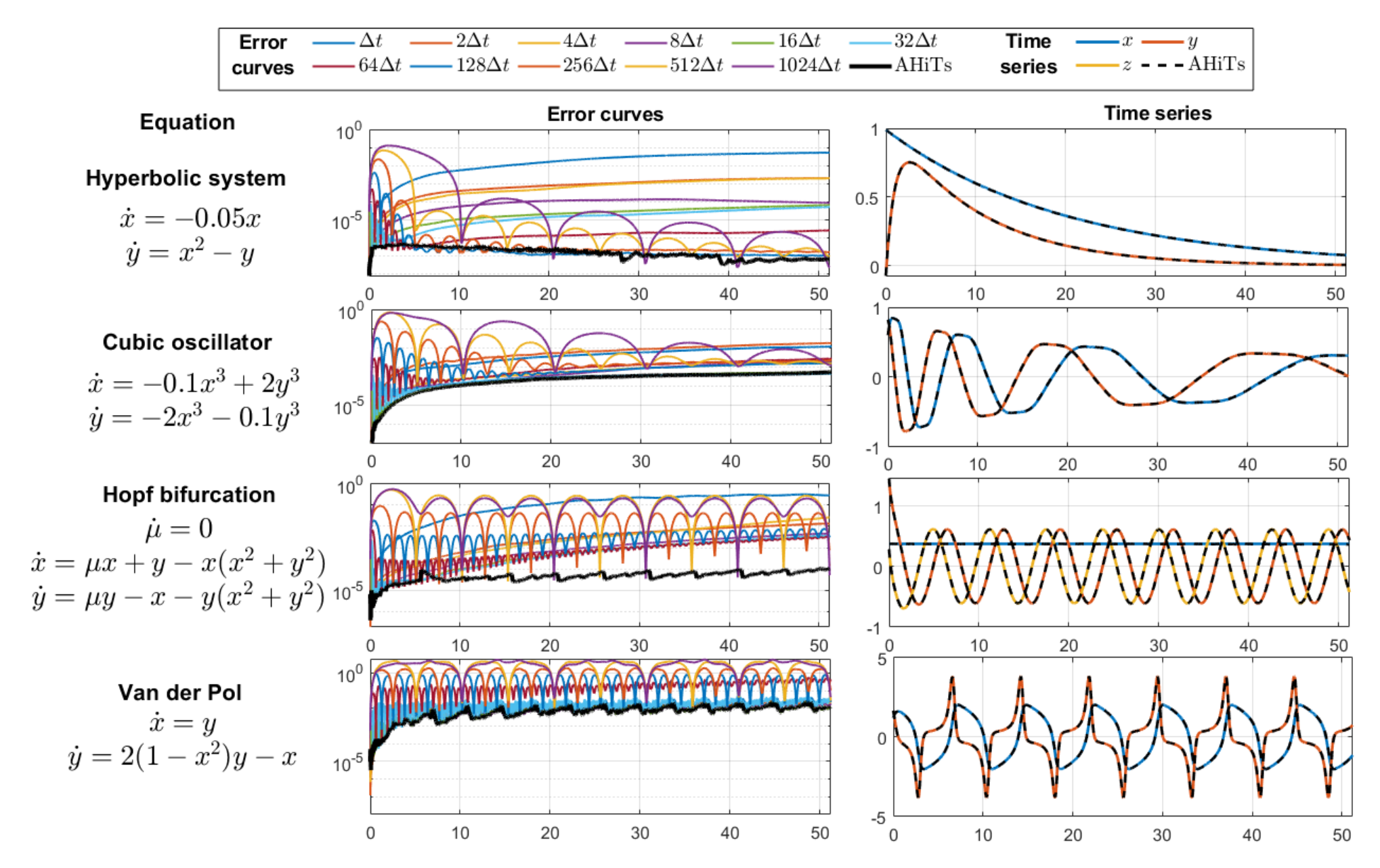}
	\caption{Performance of AHiTS vs individual NNTS models on canonical ODEs}
	\label{ahits_ode}
\end{figure}
 Figure \ref{ahits_ode} shows the performance of the proposed AHiTS scheme against all individual networks. As can be observed, the AHiTS scheme outperforms all the individual neural network responses. The error curves reveal that each network is highly accurate at the respective time steps and inaccurate elsewhere due to interpolation. Besides, 
the networks with smaller time steps yield accurate short-term predictions; however, the error quickly accumulates, as explained earlier. In contrast, networks with larger time steps have better long-term predictions but fail to capture information between steps. The proposed scheme balances these two factors and provides the best reconstruction than any individual network. This is also verified from Table \ref{errortable} wherein we report the $\mathcal{L}_{2}$ error averaged over all time steps and test trajectories. The computation time, on the other hand, for all individual networks and the proposed method is reported in Table \ref{ttode}. While it's not surprising that computation accelerates as step size grows, the proposed method is faster than the finest network across all test systems. This shows the efficiency of the AHiTS method. Finally, on comparing the AHiTS method with the multiscale HiTS method, as reported in Table \ref{tableode}, we observe the increased efficiency offered by the AHiTS method over the multiscale HiTS method. This is due to adaptive time stepping which reduces unnecessary steps while maintaining accuracy at par with the multiscale HiTS method.

\begin{table}[h!]
\begin{center}
	\caption{Relative mean square errors for predicted and exact measurements (noise-free)}
	\label{errortable}
	\begin{tabular}{@{}lcccc@{}}
		\toprule
		Systems & Hyperbolic & Cubic oscillator & Van der Pol & Hopf bifurcation \\
		\midrule
		NNTS 0 & $2.50e-2$ & $5.23e-2$ & $1.37e-2$ & $1.61e-1$ \\
		
		NNTS 1 & $10.4e-3$ & $7.77e-2$ & $1.32e-2$ & $4.79e-3$ \\
		
		NNTS 2 & $9.04e-4$ & $8.54e-4$ & $1.15e-2$ & $5.99e-3$ \\
		
		NNTS 3 & $9.10e-5$ & $7.56e-4$ & $9.23e-3$ & $1.53e-3$ \\
		
		NNTS 4 & $2.82e-5$ & $2.96e-4$ & $8.66e-3$ & $1.05e-3$ \\
		
		NNTS 5 & $1.98e-5$ & $2.90e-4$ & $1.93e-2$ & $1.10e-3$ \\
		
		NNTS 6 & $6.07e-6$ & $1.19e-3$ & $1.93e-1$ & $1.09e-3$ \\
		
		NNTS 7 & $5.79e-5$ & $1.70e-3$ & $4.09e-1$ & $3.54e-3$ \\
		
		NNTS 8 & $5.83e-4$ & $1.05e-2$ & $1.03e+0$ & $2.38e-2$ \\
		
		NNTS 9 & $3.14e-3$ & $5.09e-2$ & $2.88e+0$ & $1.51e-1$ \\
		
		NNTS 10 & $1.12e-2$ & $1.10e-1$ & $3.20e+0$ & $1.32e-1$ \\
		
		AHiTS & $1.53e-7$ & $2.72e-4$ & $7.83e-3$ & $4.99e-5$ \\
		\bottomrule
	\end{tabular}
\end{center}
\end{table}

\begin{table}[h!]
	\begin{center}
	\caption{Computation time for AHiTS and individual NNTS}
	\label{ttode}
	\begin{tabular}{@{}lcccc@{}}
		\toprule
		Systems & Hyperbolic & Cubic oscillator & Van der Pol & Hopf bifurcation \\
		\midrule
		NNTS 0 & $2.72s$ & $4.65s$ & $3.08s$ & $2.78s$ \\
		NNTS 1 & $1.66s$ & $2.14s$ & $1.44s$ & $1.59s$ \\

		NNTS 2 & $0.72 s$ & $1.08s$ & $0.66 s$ & $0.71 s$ \\

		NNTS 3 & $0.38s$ & $0.56s$ & $0.36s$ & $0.46s$ \\

		NNTS 4 & $0.20s$ & $0.29s$ & $0.19s$ & $0.24s$ \\

		NNTS 5 & $1.29s$ & $1.59s$ & $2.78s$ & $0.14s$ \\

		NNTS 6 & $0.07s$ & $0.11s$ & $0.09s$ & $0.11s$ \\

		NNTS 7 & $0.08 s$ & $0.07 s$ & $0.06 s$ & $0.08 s$ \\

		NNTS 8 & $0.05s$ & $0.05s$ & $0.05s$ & $0.07s$ \\

		NNTS 9 & $0.06s$ & $0.06s$ & $0.05s$ & $0.06s$ \\

		NNTS 10 & $0.03s$ & $0.05s$ & $0.05s$ & $0.08s$ \\

		AHiTS & $0.18s$ & $0.42s$ & $0.23s$ & $0.20s$ \\
		\bottomrule
	\end{tabular}
	\end{center}
\end{table}

\begin{table}[h!]
	\caption{Comparison of multiscale HiTS and proposed AHiTS method for benchmark ODEs}	\label{tableode}
	\centering
	\begin{tabular}{@{}lcccccc@{}}
		\toprule
		\multirow{2}{*}{Systems} & \multicolumn{3}{c}{multiscale HiTS \cite{liu2022hierarchical}} & \multicolumn{3}{c}{proposed AHiTS} \\
		\cline {2 - 7} 
		&Steps&CPU time&MSE&Steps&CPU time&MSE\\
		\midrule
		Hyperbolic & $656$  & $25.80s$ &$1.58e-7$  & $364$ & $2.88s$ &$1.53e-7$\\
	
		Cubic & $644$ & $48.87s$ & $2.85e-4$& $469$ & $6.17s$ &$2.72e-4$\\
	
		Van der Pol & $667$ & $24.65s$ &$8.26e-3$ &$639$ & $4.58s$ &$7.83e-3$\\
	
		Hopf Bifurcation & $1536$ & $26.77s$ &$5.35e-5$ &$326$ & $2.64s$ &$4.99e-5$\\
		\bottomrule
	\end{tabular}
\end{table}
\subsection{Benchmark PDEs}
Besides an efficient and accurate integration of simple nonlinear systems, we demonstrate that the AHiTS method can be used to predict the flow of multiscale PDEs accurately. For demonstrate this idea, we use the FitzHugh-Nagumo (FHN) model and the chaotic Kuramoto-Sivashinsky (KS) equation, which are discussed next.
\subsubsection{FitzHugh-Nagumo model}
The FitzHugh-Nagumo model is often used to study the behavior of excitable systems, such as neurons or cardiac cells \cite{fitzhugh1961impulses}. The coupled PDE is given as:
\begin{equation}
	\begin{split}
		\beta \frac{\partial \mathbf{u}}{\partial t}(\mathbf{x},t) &= \beta^2 \frac{\partial \mathbf{u}^{2}}{\partial \mathbf{x}^{2}}(\mathbf{x},t) + \mathbf{f}(\mathbf{u}(\mathbf{x},t))-\mathbf{v}(\mathbf{x},t)+ 0.05  \\
		\frac{\partial \mathbf{v}}{\partial t}(\mathbf{x},t) &= 0.5\mathbf{u}(\mathbf{x},t)- 2 \mathbf{v}(\mathbf{x},t)+0.05,
		\label{fhneq}
	\end{split}
\end{equation}
where, $\mathbf{f}(\mathbf{u})=\mathbf{u}(\mathbf{u}-0.1)(1-\mathbf{u})$. The variable $\mathbf{u}(\mathbf{x},t)$ represents an \textit{activator}, which drives the system towards excitation, while the variable $\mathbf{v}(\mathbf{x},t)$ represents an \textit{inhibitor}, which counteracts the effect of $\mathbf{u}(\mathbf{x},t)$ and helps to bring the system back to its resting state. The cubic nonlinearity gives rise to a nonlinear feedback mechanism that amplifies small deviation from the resting state and the bifurcation parameter $\beta=0.015$ adjusts the difference in timescales \cite{fitzhugh1961impulses,nagumo1962active}. For this case, we data was collected using the FDM scheme 
by discretizing the spatial domain $\mathbf{x} \in [0,1]$ into $100$ grid points for varying initial states of training, validation and testing trajectories. The various parameters used during training are enlisted in Appendix Table \ref{parPDE}.
 \par Figure \ref{ahits_pde} depicts the performance of proposed AHiTS method in capturing the FHN dynamics for the desired outputs $\mathbf{y}_{1}(t)$ and $\mathbf{y}_{2}(t)$ measured on the right boundary ($\mathbf{x}=1$) of the activator and inhibitor densities respectively. As can be seen, the proposed AHiTS method satisfactorily captures the slow and fast dynamics of the FHN model. Besides, AHiTS method outperforms all the individual NNTS models in terms of the accuracy as reported in Table \ref{pde1}. Upon comparison with the multiscale HiTS method, as reported in Table \ref{tablepde}, we observe that AHiTS requires almost half the number of steps than multiscale HiTS method to maintain the same level of accuracy. This demonstrates the efficiency of the proposed strategy.
\begin{figure}[h!]
	\centering
	\includegraphics[width=\textwidth]{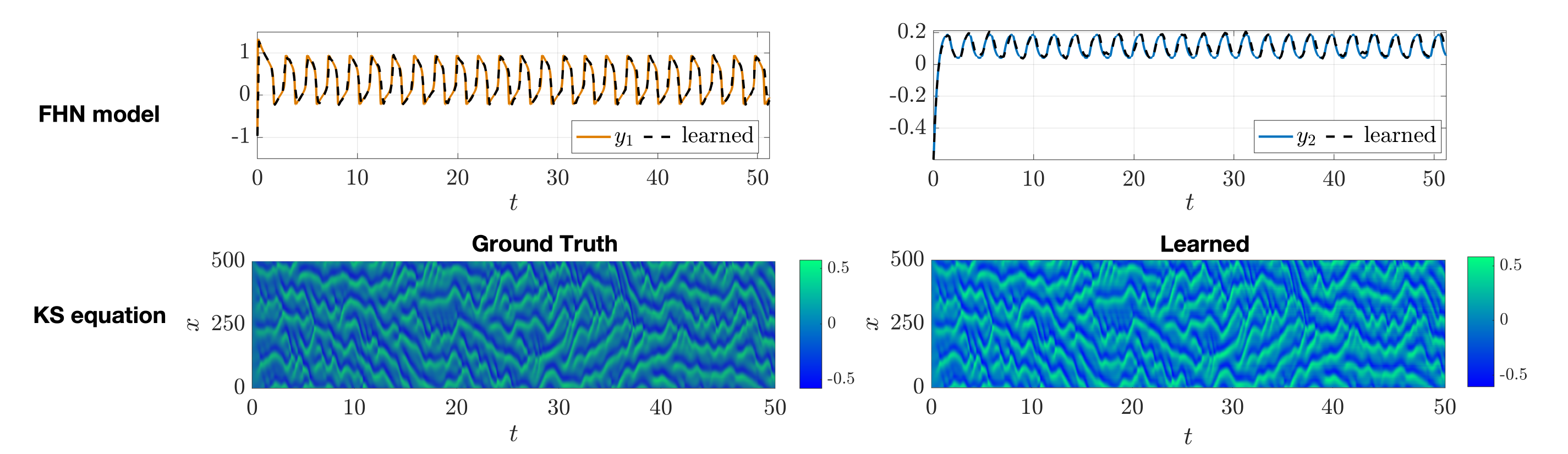}
	\caption{Performance of AHiTS on canonical PDEs}
	\label{ahits_pde}
\end{figure}
\begin{table}[h!]
	\centering
	\caption{FHN model: comparison of computation time and MSE between AHiTS and individual NNTS models}
	\label{pde1}
	\begin{tabular}{@{}lcc@{}}
		\toprule
		Systems & CPU Time & MSE \\
		\midrule
		NNTS 0 & $1.73s$ & $2.04e-1$ \\
		
		NNTS 1 & $1.36s$ & $4.08e-2$ \\
		
		NNTS 2 & $1.21s$ & $1.21e-1$ \\
		
		NNTS 3 & $1.42s$ & $8.83e-2$ \\
		
		NNTS 4 & $1.07s$ & $1.54e-1$ \\
		
		NNTS 5 & $1.32s$ & $1.32e-1$ \\
		
		NNTS 6 & $1.04s$ & $1.70e-1$ \\
		
		AHiTS & $12.56s$ & $5.20e-3$ \\
		\bottomrule
	\end{tabular}
\end{table}

\subsubsection{Kuramoto-Sivashinsky equation}
Finally, we take the one-dimensional Kuramoto-Sivashinsky (KS) equation given as:
\begin{equation}
	\mathbf{u}_t + \mathbf{u}_{xx} + \mathbf{u}_{xxxx} +\frac{1}{2} \mathbf{u}_{x}^{2}=0.
	\label{ks}
\end{equation}
The nonlinear advection term leads to the formation of coherent structures, such as fronts and shocks, which evolve on a fast timescale, whereas the diffusion term acts to smooth out the field $\mathbf{u}$ and is responsible for the long-term behavior of the system. The fourth-order dispersion term introduces an additional spatial scale and can lead to the formation of secondary structures, such as ripples and wrinkles, on a slower timescale \cite{kuramoto1978diffusion,sivashinsky1977nonlinear}.
\par 
 In this case, we considered a spatial discretization of $512$ grid points for a time span of $t\in [0,50]$. The various parameters for this case are presented in Appendix Table \ref{parPDE}. From Fig. \ref{ahits_pde}, one can visually see that the proposed AHiTS provides an accurate reconstruction of the original system dynamics, an these are confirmed by the $\mathcal{L}_{2}$ errors shown in Table \ref{pde2}. We also see that the proposed AHiTS scheme performs competitively better than the multiscale HiTS method as reported in Table \ref{tablepde}.

	\begin{table}
		\centering
			\caption{KS equation: comparison of the computation time and MSE between AHiTS scheme and individual NNTS models}
		\label{pde2}
		\begin{tabular}{@{}lcc@{}}
			\toprule
			Systems  & CPU Time & MSE \\
			\midrule
			NNTS 0 & $0.88s$ & $4.48e-1$ \\
			
			NNTS 1 & $0.41s$ & $2.22e-2$ \\
			
			NNTS 2 & $0.21s$ & $2.28e-3$ \\
		
			NNTS 3 & $0.12s$ & $1.81e-3$ \\
		
			NNTS 4 & $0.06s$ & $5.19e-2$ \\
		
			NNTS 5 & $0.04s$ & $4.03e-2$ \\
		
			NNTS 6 & $0.05s$ & $2.84e-2$ \\
		
			NNTS 7 & $0.03s$ & $4.72e-2$ \\
		
			NNTS 8 & $0.02s$ & $6.92e-2$ \\
		
			NNTS 9 & $0.01 s$ & $7.64e-2$ \\
		
			NNTS 10 & $0.03s$ & $8.12e-2$ \\
			
			AHiTS & $15.91s$ & $2.64e-5$ \\
			\bottomrule
		\end{tabular}
\end{table}

\begin{table}[h!]
	\centering
	\caption{Comparison of multiscale HiTS and the proposed AHiTS method for canonical PDEs}
	\label{tablepde}
	\begin{tabular}{@{}lcccccc@{}}
		\toprule
		\multirow{2}{*}{Systems} & \multicolumn{3}{c}{multiscale HiTS \cite{liu2022hierarchical}} & \multicolumn{3}{c}{proposed AHiTS} \\
		\cline {2 - 7}
		&Steps&CPU time&MSE&Steps&CPU time&MSE\\
		\midrule
		FHN model & $640$ & $22.94s$ & $5.70e-3$ & $336$ & $12.56s$ & $5.20e-3$ \\
		KS equation   & $4096$ & $3.50s$ & $2.64e-5$ & $4000$ & $15.91s$ & $2.64e-5$ \\
		\bottomrule
	\end{tabular}
\end{table}

\section{Discussion and Outlook}
\label{sec-discuss}
This manuscript presents a novel adaptive time stepping scheme based on training deep neural networks
for multiscale systems. Our method outperforms the fixed step size-based multiscale time stepping in terms of computational efforts. We have demonstrated that using adaptive time stepping can allow to capture slow and fast dynamics of the system effectively. We also discussed how an iterative vectorization can benefit from reducing the local error propagation of various networks. We validated the proposed method on several canonical dynamical systems and on some high-dimensional problems. 
In the following, we discuss certain aspects of the proposed scheme and then provide an outlook of this method at the end.
\subsection{Effect of noise}
\label{sec-noise}
In Tables \ref{n1},\ref{n2},\ref{n5},\ref{n10}, and \ref{n20}, we study the accuracy of various NNTS models and the robustness of proposed method against noise to depict real-world measurements. We add Gaussian random noise with varying level of variances to each component of the dynamics. The variances are set to be $1\%$, $2\%$, $5\%$, $10\%$ and $20\%$ of the variance of that component averaged over all trajectories across the data sets. We observe that NNTS models with larger temporal gaps are more resilient to noise than networks with smaller temporal gaps. This is consistent with the findings in Refs. \cite{raissi2019physics,liu2022hierarchical}. The reason is that a significant time gap between consecutive snapshots (i.e., a large $\Delta t$) allows more information to be captured, as the snapshots are more dissimilar. Conversely, if $\Delta t$ is too small, the importance of the neural network diminishes, making it impractical to train the model. However, the multiscale HiTS method and the proposed scheme remain consistently accurate than any NNTS model as they benefits from the hierarchical learning which combines the small as well as large temporal gaps.
\subsection{Choice of tolerance}
\label{tol}
The error tolerance $\epsilon$ is a key parameter that affects the proposed scheme's accuracy and efficiency. The AHiTS method adjusts the NNTS models based on the local error estimate to achieve a specified level of accuracy while minimizing the computational cost. If the error tolerance is too high, the AHiTS method may terminate before reaching the desired accuracy, resulting in inaccurate results. This would then encourage bigger time steps and be more efficient but at the expense of accuracy.
On the other hand, if the error tolerance is set too low, the AHiTS method may use more computational resources taking unnecessarily small steps, resulting in inefficient calculations. Finding the optimal error tolerance depends on several factors, such as the system's complexity, the desired level of accuracy, and the available computational resources. The choices of tolerance values used for all test cases are mentioned in Appendix Tables \ref{parode} and \ref{parPDE}. 

\subsection{Current limitations and future applications}
Although the proposed method works efficiently for small-scale and large-scale systems, it entails a substantial computational cost for training individual networks. This contrasts with the classical numerical time integrators that derive the state measurements by querying the known vector field at a few locations. However, as discussed before, these are limited by the \textit{Taylor} series constraints. Thus, an effective way to overcome this bottleneck would be to use a \textit{hybrid} approach (as suggested in Ref. \cite{liu2022hierarchical}), wherein the small time steps are recorded by the classical numerical schemes and the larger time steps are provided using a neural network, making the task inherently parallelizable. Similar to this approach, the proposed method can be combined with variable step solvers to make them computationally efficient for long-term integration. In the case of PDEs, the training costs are even higher. This is due to the large volume of data obtained during spatial discretization. A natural remedy to this would be to use an \textit{encoder-decoder} framework, wherein the state measurements will be first compressed before training the NNTS models. 
\par   
Apart from high training costs, some cautionary remarks exist when using neural networks to learn flow maps. First, it is essential to ensure that the data used to train the neural network is representative of the system's behavior over a wide range of conditions. Second, the neural network should be validated against independent data to ensure it generalizes well and does not overfit the training data. Finally, the neural network should be used with care when extrapolating beyond the range of the training data, as it may produce unreliable results.

\section*{Credit authorship contribution statement}
\textbf{Asif Hamid:} Software, Methodology, Writing- original draft,
\textbf{Danish Rafiq:} Conceptualization, Investigation, Software, Writing- original, Reviewing and Editing,
\textbf{Shahkar A. Nahvi:}  Supervision, Reviewing,
\textbf{Mohammad A. Bazaz:} Supervision, Reviewing
\section*{Declaration of competing interest} The authors declare no potential financial or non-financial competing interests
\section*{Acknowledgments}
Author I would like to acknowledge the Doctoral fellowship from the Ministry of Education (MoE), New Delhi, India via Grant No. IUST0119013135, and Author II acknowledges the financial assistance from the Science and Engineering Research Board (SERB), a statutory body of the Department of Science and Technology (DST), Government of India via Grant No. PDF/2022/002081.
\section*{Code availability}
The source code is available at: \url{https://github.com/DanishRaf32/Adaptive-HiTs}
\bibliography{mybib}
\hspace{10pt}

\begin{appendices}
\section{Algorithms}\label{sec-algo}
\begin{algorithm}[H]
	\caption{: Adaptive time step estimation}
	\label{alglo1}
\textbf{Input:}	Set of neural-network time steppers $\mathcal{N}_d (.: \Delta t_{d})$, tolerance $\epsilon$, unit step size $\Delta t$, step sizes for individual NNTS:  $stepsizes=\Delta t_{d}=2^{d}\Delta t~ (d=0,..,m)$, final time =$t_{f}$, initial conditions $\mathbf{x}_0 \in \mathbb{R}^{\boldsymbol{\mu} \times 1 \times n}$ from validation data set\\ 
	\textbf{Output:} Adaptive time steps: $steps$ \\[-0.4cm]
	\begin{algorithmic}[1]
		\State $ \mathcal{N}_d=\textit{sort}(\mathcal{N}_d$)   \hfill $\triangleright$ sort NNTS in descending order.
		\State Initialize $y_{prev}=\mathbf{x}_0, k=1, steps=[~]$  
		\While{$k \leqslant {t_f}/\Delta t$} 
		\For {$i = 1,..,m$} 
			\State $y_{next}= \mathcal{N}_i (y_{prev})$
			\State $mse=\vert \vert y_{next} -y_{prev} \vert \vert^{2} $ \hfill $\triangleright$ calculate the evolution of states
			\If{$mse < \epsilon$}
				\State $y_{next}= y_{prev}$
				\State $k=sum(stepsizes(i))$         \hfill  $\triangleright$ add current time step
				\State $steps=append(stepsizes(i))$   \hfill  $\triangleright$ append current time step
				\State break
			\EndIf
			\If{$i==m$}	 \hfill  $\triangleright$ check if no NNTS is selected
				\State $y_{next}= y_{prev}$
				\State $k=sum(stepsizes(i))$ 
				\State $steps=append(stepsizes(i))$ 
			\EndIf
		\EndFor
		\EndWhile
	\end{algorithmic}
\end{algorithm}

\begin{algorithm}[H]
	\caption{: Iterative vectorized computation}
	\textbf{Input:}	Adaptive steps: $steps$, Neural network time steppers $\mathcal{N}_{d}(.; \Delta t_{d})$, step sizes: $stepsizes=\Delta t_{d}=2^{d}\Delta t$, validation data $\mathbf{X} \in \mathbb{R}^{\boldsymbol{\mu} \times t_f \times n}$ \\ 
	\textbf{Output:}  shortlisted models  $\mathcal{N}_v$\\[-0.4cm]
	\label{alglo2}
	\begin{algorithmic}[1]
		\State  $count=0, curr\_steps=0, k=0, j=1, iid=[~]$ \hfill $\triangleright$ initialize variables
		\For{$i = 1,..,len(steps)$} 	
	     	\State $x=steps(i)$     \hfill $\triangleright$  store current step
			\If{$x==steps(i+1)$}	\hfill 	$\triangleright$ check if a step gets repeated
				\State $count= count+1, j=j+1$	  
			\Else
					\If{$count==1$} \hfill $\triangleright$ check if a step is used only once
				\State $\mathcal{N}_v=models[index(stepsizes(x))]$
				\Else
				\State $\mathcal{N}_v=models_{init}[:end_{idx}]$
				\State  $k=k+1, count=1$
				\EndIf
				
				\State $n_{steps}={count \ast steps(i-1)}-1$  
				\State $curr\_step=n_{steps}+curr\_step$
				\State $iid=append(curr\_step)$
				\State $init=index(stepsizes(i-1))$
				\State $fin=index(\arg \min (stepsizes-n_{steps}))$
				\State $models_{init}=models(init:fin)$ \hfill 	$\triangleright$ store models for one window
				\State $ic=iid(k)+1$
				\State $best_{mse}=1e+5$
				\For{$j=1,2, \dots, len(models_{init})$}
					\State $y_{pred}=vectorized(models_{init}[:len(models_{init})-j])$	\hfill  $\triangleright$   Ref.  \cite{liu2022hierarchical}
					\State $mse=\mathbf{X}^{val}(ic:n_{steps}+1)-y_{pred}$
					\If{$mse \leqslant best_{mse} $}
						\State $end_{idx}=len(models_{init})-j$
						\State $best_{mse}=mse$
					\EndIf
		   	     \EndFor

						\State $\mathcal{N}_v=append(\mathcal{N}_v)$ 
		   \EndIf
		\EndFor
	\end{algorithmic}
\end{algorithm}

\begin{sidewaystable}[h!] 
	\section{Parameters}	
	\label{sec-params}
	\caption{Parameters for various network architectures used for ODEs}
	\label{parode}
	\begin{tabular*}{\textwidth}{@{}lccccc@{}}
		\toprule 
		Systems &\begin{tabular}{@{}c@{}}
			Samples\\ (train/\\validation/\\test)
		\end{tabular} &\begin{tabular}{@{}c@{}}
		Sampled region\\ $\mathcal{D}$ in state space
	\end{tabular} & \begin{tabular}{@{}c@{}}
	Network\\ architectures
\end{tabular} & \begin{tabular}{@{}c@{}}
Cross-validated\\ NNTS for HiTS 
\end{tabular} &\begin{tabular}{@{}c@{}}
Tolerance for \\ AHiTS
\end{tabular} 
		\\
		\midrule
		Hyperbolic &  $1600/320/320$ & $[-1,1]^{2}$ &$[2,128,128,128,2]$  & $3-5$ & $1e-5$ \\
	
		Cubic oscillator  & $3200/320/320$ & $[-1,1]^{2}$ & $[2,256,256,256,2]$ & $3-6$ & $5e-4$ \\
	
		Van der Pol &$3200/320/320$  & $[-2,2] \times[-4,4]$ & $[2,512,512,512,2]$ & $3-5$ & $8e-2$ \\
	
		Hopf bifurcation & $3200/320/320$  & $[-0.2,0.6] \times[-1,2] \times[-1,1]$ & $[3,128,128,,128,3]$ & $2-10$ & $5e-3$ \\
	\bottomrule
	\end{tabular*}
\vspace{0.5cm}
	\caption{Parameters for various network architectures used for PDEs}
	\label{parPDE}
	\begin{tabular}{@{}lcccc@{}}
		\toprule
		Systems & Data size & \begin{tabular}{@{}c@{}}
			Network\\ architectures
		\end{tabular} & \begin{tabular}{@{}c@{}}
			Cross-validated\\ NNTS for HiTS 
		\end{tabular} &\begin{tabular}{@{}c@{}}
			Tolerance for \\ AHiTS \end{tabular} \\
		\midrule
		KS equation &$1 \times 4001 \times 512 $&\begin{tabular}{@{}c@{}}
			 $[512,2048,512], [512,1024,512], [512,512,512],$ \\ $ [512,256,512], [512,128,512]$ 
		\end{tabular}  & $0-10$ & $1e-5$ \\
	
	FHN model & $3300 \times 5121 \times 100$ & \begin{tabular}{@{}c@{}}$[100,100,512,1024,2048,1024,512,100,100]$\end{tabular} & $0-5$ & $2e-2$ \\
		\bottomrule
	\end{tabular}
\end{sidewaystable}

\begin{table}[h!] 
	\section{Noisy measurements}
	\label{sec-noisy}
	\caption{Comparison of the error response for $1 \%$ noise}
	\label{n1}
	\centering
	\begin{tabular}{@{}lcccc@{}}
		\toprule
		Systems & Hyperbolic & Cubic oscillator & Van der Pol & Hopf bifurcation \\
		\midrule
		NNTS 0 & $2.89e-2$ & $1.12e-1$ & $4.15e+0$ & $2.79e+2$ \\
	
		NNTS 1 & $3.40e-3$ & $1.81e-1$ & $2.25e+0$ & $2.11e-2$ \\
		
		NNTS 2 & $8.36e-4$ & $2.23e-2$ & $2.73e+0$ & $3.66e-3$ \\
		
		NNTS 3 & $1.82e-4$ & $2.36e-2$ & $1.48e+0$ & $5.64e-3$ \\
		
		NNTS 4 & $1.03e-4$ & $5.70e-3$ & $5.15e-1$ & $1.21e-3$ \\
		
		NNTS 5 & $1.47e-5$ & $1.03e-2$ & $5.15e-1$ & $3.74e-2$ \\
		
		NNTS 6 & $1.19e-5$ & $4.84e-2$ & $1.38e-1$ & $6.53e-4$ \\
	
		NNTS 7 & $6.17e-5$ & $2.53e-3$ & $4.16e-1$ & $2.17e-3$ \\
		
		NNTS 8 & $5.81e-4$ & $1.11e-2$ & $1.02e+0$ & $2.29e-2$ \\
	
		NNTS 9 & $3.42e-3$ & $5.35e-2$ & $2.85e+0$ & $1.51e-1$ \\
		
		NNTS 10 & $1.12e-2$ & $1.02e-1$ & $3.21e+0$ & $1.31e-1$ \\
		
		AHiTS & $1.86e-6$ & $1.46e-3$ & $2.47e-2$ & $5.72e-5$ \\
		\bottomrule
	\end{tabular}

\end{table}


\begin{table}[h!]
	\caption{Comparison of the error response for $2 \%$ noise}
	\label{n2}
	\centering
	\begin{tabular}{@{}lcccc@{}}
	\toprule
	Systems & Hyperbolic & Cubic oscillator & Van der Pol & Hopf bifurcation \\
	\midrule
		NNTS 0 & $1.75e-1$ & $1.11e-1$ & $2.32e+0$ & $2.81e-1$ \\
	
		NNTS 1 & $1.72e-2$ & $8.35e-2$ & $5.22e-1$ & $6.53e-1$ \\
	
		NNTS 2 & $6.43e-3$ & $9.31e-2$ & $2.89e-1$ & $1.32e-2$ \\
	
		NNTS 3 & $9.28e-4$ & $5.86e-2$ & $9.14e-1$ & $1.79e-2$ \\
	
		NNTS 4 & $2.89e-4$ & $5.56e-3$ & $1.50e+0$ & $2.37e-3$ \\
	
		NNTS 5 & $1.99e-4$ & $3.56e-2$ & $7.52e-1$ & $8.35e-4$ \\
	
		NNTS 6 & $2.24e-5$ & $9.51e-3$ & $1.59e-1$ & $6.12e-4$ \\
		
		NNTS 7 & $7.44e-5$ & $3.25e-3$ & $4.58e-1$ & $2.57e-3$ \\
		
		NNTS 8 & $5.48e-4$ & $3.31e-2$ & $1.04e+0$ & $2.33e-2$ \\
	
		NNTS 9 & $3.21e-3$ & $5.29e-2$ & $2.82e+0$ & $1.51e-1$ \\
	
		NNTS 10 & $1.05e-2$ & $1.03e-1$ & $3.23e+0$ & $1.32e-1$ \\
	
		AHiTS & $7.25e-6$ & $2.94e-3$ & $8.76e-2$ & $9.26e-5$ \\
		\bottomrule
	\end{tabular}

\end{table}


\begin{table}[h!]
	\caption{Comparison of the error response for $5 \%$ noise}
	\label{n5}
	\centering
	\begin{tabular}{@{}lcccc@{}}
	\toprule
	Systems & Hyperbolic & Cubic oscillator & Van der Pol & Hopf bifurcation \\
	\midrule
		NNTS 0 & $4.34e-2$ & $2.08e-1$ & $3.20e+0$ & $8.78e-2$ \\
	
		NNTS 1 & $8.31e-2$ & $1.37e-1$ & $3.47e+0$ & $4.28e-1$ \\
	
		NNTS 2 & $7.83e-3$ & $4.75e-2$ & $3.15e+0$ & $4.37e-2$ \\

		NNTS 3 & $2.40e-3$ & $4.30e-2$ & $1.20e+0$ & $6.99e-3$ \\
	
		NNTS 4 & $2.10e-3$ & $2.12e-2$ & $1.89e+0$ & $4.71e-3$ \\
		
		NNTS 5 & $2.47e-4$ & $1.09e-2$ & $3.06e+0$ & $1.29e-3$ \\
		
		NNTS 6 & $5.38e-5$ & $4.17e-2$ & $6.24e-1$ & $9.62e-4$ \\
	
		NNTS 7 & $7.11e-4$ & $1.32e-2$ & $4.93e-1$ & $4.54e-3$ \\
	
		NNTS 8 & $5.77e-4$ & $2.56e-2$ & $1.06e+0$ & $2.37e-2$ \\

		NNTS 9 & $3.24e-3$ & $5.85e-2$ & $2.85e+0$ & $1.51e-1$ \\
	
		NNTS 10 & $1.05e-2$ & $1.06e-1$ & $3.17e+0$ & $1.32e-1$ \\
		
		AHiTS & $3.96e- 5$ & $8.90e-3$ &$1.84e-1$ & $3.41e-4$ \\
		\bottomrule
	\end{tabular}
\end{table}

\begin{table}[t]
	\caption{Comparison of the error response for $10 \%$ noise}
	\label{n10}
	\centering
		\begin{tabular}{@{}lcccc@{}}
		\toprule
		Systems & Hyperbolic & Cubic oscillator & Van der Pol & Hopf bifurcation \\
		\midrule
		NNTS 0 & $4.95e-1$ & $4.45e-1$ & $3.62e+0$ & $2.90e-1$ \\
	
		NNTS 1 & $1.08e-1$ & $2.10e-1$ & $3.84e+0$ & $2.62e-1$ \\
	
		NNTS 2 & $7.23e-2$ & $1.82e-1$ & $5.22e+0$ & $1.35e-1$ \\

		NNTS 3 & $1.46e-3$ & $7.12e-2$ & $4.66e+0$ & $1.16e-2$ \\
	
		NNTS 4 & $9.67e-4$ & $4.38e-2$ & $4.16e+0$ & $2.68e-2$ \\
	
		NNTS 5 & $1.52e-3$ & $4.19e-2$ & $4.33e+0$ & $8.18e-3$ \\

		NNTS 6 & $4.89e-4$ & $6.69e-2$ & $7.96e+0$ & $3.56e-3$ \\
	
		NNTS 7 & $5.10e-4$ & $2.80e-2$ & $7.31e-1$ & $3.12e-3$ \\
	
		NNTS 8 & $7.30e-4$ & $4.99e-2$ & $1.7e+0$ & $2.39e-2$ \\
	
		NNTS 9 & $3.29e-3$ & $7.59e-2$ & $2.85e+0$ & $1.52e-1$ \\
	
		NNTS 10 & $1.06e-2$ & $1.14e-1$ & $3.15e+0$ & $1.32e-1$ \\
	
		AHiTS & $1.24e-4$ & $2.25e-2$ & $4.02e-1$ & $1.14e-3$ \\
		\bottomrule
	\end{tabular}
\vspace{1.5cm}

	\caption{Comparison of the error response for $20 \%$ noise}
	\label{n20}
	\centering
		\begin{tabular}{@{}lcccc@{}}
		\toprule
		Systems & Hyperbolic & Cubic oscillator & Van der Pol & Hopf bifurcation \\
		\midrule
		NNTS 0 & $1.68e-1$ & $4.73e-1$ & $3.31e+0$ & $1.97e-1$ \\
		
		NNTS 1 & $1.48e-1$ & $2.04e-1$ & $2.80e+0$ & $1.83e-1$ \\
	
		NNTS 2 & $1.50e-1$ & $1.60e-1$ & $2.95e+0$ & $2.17e-1$ \\
	
		NNTS 3 & $5.97e-2$ & $2.74e-1$ & $4.29e+0$ & $5.61e-2$ \\
	
		NNTS 4 & $7.05e-3$ & $1.21e-1$ & $2.56e+0$ & $4.22e-2$ \\
		
		NNTS 5 & $1.24e-3$ & $9.81e-2$ & $3.99e+0$ & $2.73e-2$ \\
		
		NNTS 6 & $7.73e-4$ & $1.02e-1$ & $3.99e+0$ & $1.23e-2$ \\
	
		NNTS 7 & $8.08e-4$ & $9.99e-2$ & $2.39e+0$ & $7.79e-3$ \\
	
		NNTS 8 & $11.3e-3$ & $9.35e-2$ & $1.21e+0$ & $2.80e-2$ \\
	
		NNTS 9 & $3.78e-3$ & $1.11e-1$ & $2.77e+0$ & $1.52e-1$ \\
	
		NNTS 10 & $1.11e-2$ & $1.33e-1$ & $3.04e+0$ & $1.34e-1$ \\
	
		AHiTS & $3.21e-4$ & $5.87e-2$ & $8.13e-1$ & $4.21e-3$ \\
		\bottomrule
	\end{tabular}
\end{table}

\end{appendices}
\end{document}